\documentclass[12pt]{article}

\usepackage{amssymb}

\makeatletter
 
 \@addtoreset{equation}{section}
\makeatother

\newcommand{\ul}{\underline}

\renewcommand{\pmod}[1]{%
\ (\mbox{\rm mod}\ #1)}

\newcommand{\prf}{\noindent{\bf Proof. }}

\newcommand{\rem}{\noindent{\bf Remark. }}

\newcommand{\ff}[1]{%
    {\bf F}_{#1}}

\newcommand{\qed}{\hbox{\rule[-2pt]{5pt}{11pt}}}

\newtheorem{dfn}{Definition}[section]
\newtheorem{thm}[dfn]{Theorem}
\newtheorem{prop}[dfn]{Proposition}

\newtheorem{lem}[dfn]{Lemma}
\newtheorem{exam}[dfn]{Example}

\newtheorem{prob}[dfn]{Problem}

\begin{document}
\title{An abundance of invariant polynomials satisfying the Riemann hypothesis}
\author{Koji Chinen\footnotemark[1]}
\date{29/4/2007}
\maketitle


\begin{abstract}
In 1999, Iwan Duursma defined the zeta function for a linear code as a  generating function of its Hamming weight enumerator. It can also be defined for other homogeneous polynomials not corresponding to existing codes. If the homogeneous polynomial is invariant under the MacWilliams transform, then its zeta function satisfies a functional equation and we can formulate an analogue of the Riemann hypothesis. As far as existing codes are concerned, the Riemann hypothesis is believed to be closely related to the extremal property. 

In this article, we show there are abundant polynomials invariant by the MacWilliams transform which satisfy the Riemann hypothesis. The proof is carried out by explicit construction of such polynomials. To prove the Riemann hypothesis for a certain class of invariant polynomials, we establish an analogue of the Enestr\"om-Kakeya theorem. 
\end{abstract}
\footnotetext[1]{Department of Mathematics, School of Science and Engineering, Kinki University. 3-4-1, Kowakae, Higashi-Osaka, 577-8502 Japan. E-mail: chinen@math.kindai.ac.jp}


\section{Introduction}

Let $p$ be a prime, $q=p^r$ for some positive integer $r$ and we denote the finite field with $q$ elements by $\ff q$. Let $C$ be an $[n,k,d]$-code over $\ff q$ with the Hamming weight enumerator $W_C(x,y)$. Duursma \cite{Du1} defined the zeta function for $C$ as a generating function of $W_C(x,y)$. Then the author \cite{Ch2} considered the case of so-called ``formal weight enumerators'', noticing that Duursma's definition can be extended for other homogeneous polynomials than the weight enumerators of actual codes. Taking these into account, we start from the following definition: 
\begin{dfn}\label{dfn:zeta}
For any $q\in{\bf N}$ ($q\geq2$) and any homogeneous polynomial of the form 
\begin{equation}\label{eq:homogen}
W(x,y)=x^n+\sum_{i=d}^n A_i x^{n-i}y^i \quad (A_i\in {\bf C},\ A_d\ne0)
\end{equation}
there exists a unique polynomial $P(T)\in{\bf C}[T]$ of degree at most $n-d$ such that
\begin{equation}\label{eq:zeta_duursma}
\frac{P(T)}{(1-T)(1-qT)}(y(1-T)+xT)^n=\cdots +\frac{W(x,y)-x^n}{q-1}T^{n-d}+ \cdots.
\end{equation}
We call $P(T)$ and $Z(T)=P(T)/(1-T)(1-qT)$ the zeta polynomial and the zeta function of $W(x,y)$, respectively. 
\end{dfn}
For the proof of existence and uniqueness of $P(T)$, see Appendix. If $W(x,y)=W_C(x,y)$ for some linear code $C$, then we take $q$ in the above definition as $\sharp \ff q$, but if $W(x,y)$ is not related to an existing code, then $q$ must be chosen suitably according to what meaning $W(x,y)$ has. 

In the case $W(x,y)=W_C(x,y)$, the zeta polynomial $P(T)$ for $W_C(x,y)$ is of particular interest when $C$ is self-dual, because it has the functional equation 
\begin{equation}\label{eq:func_eq}
P(T)=P\Bigl(\frac{1}{qT}\Bigr)q^g T^{2g}
\end{equation}
($g=n/2+1-d$, see \cite[p.59]{Du2}), which is a result of the fact that $W_C(x,y)$ is invariant by the MacWilliams transform
\begin{equation}\label{macwilliams}
\sigma_q:=\frac{1}{\sqrt q}
\left(\begin{array}{rr} 1 & q-1 \\ 1 & -1 \end{array}\right),
\end{equation}
where we define $f^\sigma(x,y)=f(ax+by,cx+dy)$ for $f(x,y)\in {\bf C}[x,y]$ and a linear transformation $\sigma=\left(\begin{array}{rr} a & b \\ c & d \end{array}\right)$. 

The functional equation (\ref{eq:func_eq}) is the same as that of zeta functions of algebraic curves, so we can formulate the Riemann hypothesis (see Duursma \cite[Definition 4.1]{Du3}). Even if $W(x,y)$ does not correspond to an actual code, we can formulate the Riemann hypothesis in the same way provided that $W^{\sigma_q}(x,y)=W(x,y)$ because it is this invariance that yields (\ref{eq:func_eq}): 
\begin{dfn}\label{dfn:RH}
The code $C$ (or the invariant polynomial $W(x,y)$) satisfies the Riemann hypothesis if all the zeros of $P(T)$ have the same absolute value $1/\sqrt q$. 
\end{dfn}
Duursma deduces various interesting properties of $P(T)$ and discusses their possible applications to the coding theory (see \cite{Du2, Du3, Du4}).

Finding an equivalent condition for the Riemann hypothesis above seems still an open problem, but Duursma asks the following (\cite[Open Problem 4.2]{Du3}): 
\begin{prob}\label{prob:RH}
Prove or disprove that all extremal weight enumerators satisfy the Riemann hypothesis. 
\end{prob}
A self-dual code $C$ is called extremal if it has the largest possible minimum distance (see Pless \cite[p.139]{Pl}). There are 4 well-known sequences of extremal self-dual codes (Types I, II, III and IV, see Conway-Sloane \cite{CoSl}). The extremal code is also characterized by its weight enumerator $W_C(x,y)$: the code $C$ is extremal if $d$ of $W_C(x,y)$ in (\ref{eq:homogen}) is the largest among all the self-dual weight enumerators of degree $n$ over $\ff q$. Using this, the extremal property is straightfowardly extended to the case of some more general invariant polynomials. Duursma proved that all extremal Type IV codes satisfied the Riemann hypothesis (\cite{Du4}). Thus, as far as the existing codes are concerned, we may expect that the Riemann hypothesis reflects one of the abilities of the code, the extremal property. 

In \cite{Ch2}, the author extended the consideration to the case of the formal weight enumerators. A formal weight enumerator $W(x,y)$ resembles the weight enumerator of a Type II code, but is distinguished from it by the property $W^{\sigma_2}(x,y)=-W(x,y)$ (see \cite[Definition 1.4]{Ch2}). The zeta polynomial $P(T)$ of $W(x,y)$ satisfies $P(T)=-P(1/2T)2^gT^{2g}$ ($g=n/2+1-d$) and we can formulate the Riemann hypothesis in the same way as in Definition \ref{dfn:RH}, setting $q=2$. In \cite[Section 3]{Ch2}, we observed that the extremal property might yield the Riemann hypothesis also in the case of the formal weight enumerators. 

The purpose of the present article is to extend the consideration to all the polynomials which are invariant by the MacWilliams transform $\sigma_q$. Such polynomials form an invariant polynomial ring 
\begin{equation}\label{eq:ipr_macwil}
{\bf C}[x,y]^{G_q}={\bf C}[x+(\sqrt{q}-1)y, y(x-y)]
\end{equation}
where $G_q=\langle \sigma_q \rangle$ (see MacWilliams-Sloane \cite[p.605, Theorem 5]{MaSl}). As a problem of invariant polynomials, we can remove the structure of linear codes and allow $q$ to be any positive integer such that $q\geq 2$. We try to find as many polynomials as possible in ${\bf C}[x,y]^{G_q}$ which satisfy the Riemann hypothesis. The results imply that the Riemann hypothesis is not always relevant to the extremal property in the ring ${\bf C}[x,y]^{G_q}$. The first result is the following: 
\begin{thm}\label{thm:RH_MDS}
For any $q\geq 2$ and any $n$, $d$ such that $2\leq d \leq \frac{n+1}{2}$, there exists a $\sigma_q$-invariant polynomial of the form (\ref{eq:homogen}) which satisfies the Riemann hypothesis. 
\end{thm}
Note that the restriction $2\leq d$ comes from the original Duursma theory. We also note that the number $d$ in (\ref{eq:homogen}) must satisfy $d \leq \frac{n}{2}+1$ in ${\bf C}[x,y]^{G_q}$. In cases where equality holds, the polynomial becomes MDS and the zeta polynomial is a constant (see Section 3). Thus by the condition $2\leq d \leq \frac{n+1}{2}$, almost all possible pairs of $n$ and $d$ are covered and it shows that polynomials satisfying the Riemann hypothesis are widely and abundantly distributed in ${\bf C}[x,y]^{G_q}$. 

The notion of extremal polynomial is defined only in terms of $n$ and $d$, but Theorem \ref{thm:RH_MDS} implies that, at least in the ring ${\bf C}[x,y]^{G_q}$, the condition for the Riemann hypothesis is not determined by $n$ and $d$ only. Thus Theorem \ref{thm:RH_MDS} shows us another aspect of the zeta functions for invariant polynomials. 

Theorem \ref{thm:RH_MDS} is proved by explicit construction of the invariant polynomials with the desired property. This is done by using the weight enumerators of codes which are not self-dual. If $C$ is not self-dual, its weight enumarator $W_C(x,y)$ does not satisfy $W_C^{\sigma_q}(x,y)=W_C(x,y)$, but combining $W_C(x,y)$ and $W_{C^\perp}(x,y)$, we can easily get an invariant expression $\tilde W_C(x,y)$ and its zeta polynomial $\tilde P_C(T)$ (see Section 2). Theorem \ref{thm:RH_MDS} is the result of the case where $C$ is an MDS code (see Section 3). 

Such a way of constructing invariant polynomials can be applied to any linear code which is not self-dual and leads to further exploration. The rest of the paper is devoted to the analysis of two other special classes of codes, the general Hamming codes and the Golay codes (not self-dual). These codes, along with certain MDS codes form an important class of good codes, the perfect codes (see Pless \cite[p.21]{Pl}): 
\begin{dfn}\label{dfn:perfect}
A code $C\subset {\ff q}^n$ of minimum distance $d$ is called perfect if all the vectors in ${\ff q}^n$ are contained in a ball of radius $[(d-1)/2]$ about the codewords, where $[x]$ means the largest integer not greater than $x$. 
\end{dfn}
The nontrivial linear perfect codes are completely determined (\cite[Section 2.2]{Pl} or \cite[Section 6.10]{MaSl}): 
\begin{itemize}
\item[(i)] The general Hamming $[(q^r-1)/(q-1)=n, n-r, 3]$ codes over $\ff q$, 
\item[(ii)] The binary $[23,12,7]$ and the ternary $[11,6,5]$ Golay codes. 
\end{itemize}
We also have trivial perfect codes: the whole space and a binary repetition code of odd length. The latter has the parameter $[n,1,n]$, being MDS and dealt with in Theorem \ref{thm:RH_MDS}. As to the general Hamming codes, they become MDS when $r=2$, so it follows that this case is also treated in Theorem \ref{thm:RH_MDS}. 

We can find again infinitely many polynomials in ${\bf C}[x,y]^{G_q}$ satisfying the Riemann hypothesis by constructing $\tilde W_C(x,y)$ from the above class of codes: 
\begin{thm}\label{thm:RH_hamming}
Let $C={\rm Ham}(r,q)$ be the Hamming $[(q^r-1)/(q-1)=n, n-r, 3]$ code over $\ff q$. If $r\geq 3$ and $q\geq 4$, then the invariant polynomial $\tilde W_C(x,y)$ in ${\bf C}[x,y]^{G_q}$ satisfies the Riemann hypothesis.
\end{thm}
To prove Theorem \ref{thm:RH_hamming}, we deduce a certain function theoretical result concerning the distribution of the zeros of a self-reciprocal polynomial: 
\begin{thm}\label{thm:ext_enest_kake}
If $f(T)=a_0+a_1 T+ \cdots +a_k T^k + a_k T^{m-k}+a_{k-1} T^{m-k+1} + \cdots +a_0 T^m$ ($m>2k$) satisfies $a_0>a_1> \cdots > a_k>0$, then all the roots of $f(T)$ lie on the unit circle. 
\end{thm}
This is, so to speak, a self-reciprocal analogue of the famous Enestr\"om-Kakeya theorem (see Theorem \ref{thm:enest_kake}). Because of the technical difficulties, Theorem \ref{thm:RH_hamming} remains unproved when $q=2,3$ and $r\geq 3$, but numerical experiments imply that the Riemann hypothesis seems to be true in these cases. 

For the Golay codes, we have the following: 
\begin{thm}\label{thm:RH_golay}
Let $C$ be the binary $[23,12,7]$ or the ternary $[11,6,5]$ Golay code.  Then the invariant polynomial $\tilde W_C(x,y)$ satisfies the Riemann hypothesis.
\end{thm}

Thus except for the binary and ternary general Hamming codes, we can prove that the invariant polynomials $\tilde W_C(x,y)$ from the perfect codes satisfy the Riemann hypothesis.  

The rest of the article is organized as follows. In Section 2, we construct an invariant polynomial $\tilde W_C(x,y)$ from the weight enumerator of a  code $C$ (which is not always self-dual) and give an explicit form of its zeta polynomial $\tilde P_C(T)$. In Section 3, we apply the results in Section 2 to the MDS code and prove Theorem \ref{thm:RH_MDS}. In Section 4, we determine the zeta polynomial $\tilde P_C(T)$ when $C={\rm Ham}(r,q)$, the general Hamming code when $r\geq 3$ and $q\geq 2$. Section 5 is devoted to an analogue of the Enestr\"om-Kakeya theorem. Here we use several results of the classical function theory. Using it, we prove the Riemann hypothesis for $\tilde W_{{\rm Ham}(r,q)}(x,y)$ ($r\geq 3$, $q\geq 4$) in Section 6. In Section 7, we consider the case of the Golay codes, and prove Theorem \ref{thm:RH_golay} by a different method to Theorem \ref{thm:RH_hamming}. 

We have been interested in the extremal property of the weight enumerators when considering the Riemann hypothesis in the context of existing self-dual codes or a little larger class of invariant polynomials which have some connections to the coding theory, that is, the formal weight enumerators. But the results in this article show that it is not always the extremal property that yields the Riemann hypothesis in the ``largest'' ring ${\bf C}[x,y]^{G_q}$. We can observe rather pathological phenomena there. We are now in a position to seek some new structures which are larger than existing codes (but smaller than ${\bf C}[x,y]^{G_q}$), in which the Riemann hypothesis indicates some distinguished properties of invariant polynomials. 

\medskip
\noindent{\it Acknowledgment.} The author would like to express his sincere gratitude to Professor Leo Murata for an abundance of valuable advice and discussion. 
\section{Invariant polynomials and their zeta functions from arbitrary  linear codes}

Let $C$ be a linear $[n,k,d]$ code over $\ff q$ and $W_C(x,y)$ be its Hamming weight enumerator. Suppose the dual code $C^\perp$ has the parameter $[n, n-k, d^\perp]$ and we assume $d, d^\perp \geq 2$. Combining $W_C(x,y)$ and the dual weight enumerator $W_{C^\perp}(x,y)$, we can easily obtain an invariant expression $\tilde W_C(x,y)$: 
\begin{prop}\label{prop:tilde_W}
Let
\begin{equation}\label{eq:tilde_W}
\tilde W_C(x,y) := \frac{1}{1+q^{k-n/2}}\{W_C(x,y)+q^{k-n/2}W_{C^\perp}(x,y)\}.
\end{equation}
Then we have $\tilde W_C^{\sigma_q}(x,y)=\tilde W_C(x,y)$, i.e. $\tilde W_C(x,y) \in {\bf C}[x,y]^{G_q}$. 
\end{prop}
\prf The proof is evident from the MacWilliams identity
$$W_C^{\sigma_q}(x,y)=q^{k-n/2}W_{C^\perp}(x,y)\quad 
{\rm or}\quad W_{C^\perp}^{\sigma_q}(x,y)=q^{n/2-k}W_C(x,y)$$
(see \cite[p.146, Theorem 13]{MaSl}). \qed

\medskip
Now we deduce the explicit form of the zeta polynomial $\tilde P_C(T)$ of $\tilde W_C(x,y)$. Let $P_C(T)$ and $P_{C^\perp}(T)$ be the zeta polynomials of $W_C(x,y)$ and $W_{C^\perp}(x,y)$, respectively. Our goal in this section is to prove the following: 
\begin{thm}\label{thm:tilde_P}
The zeta polynomial $\tilde P_C(T)$ of $\tilde W_C(x,y)$ is given by 
\begin{equation}\label{eq:tilde_P}
\tilde P_C(T)=\frac{T^{\max(0,d-d^\perp)}}{1+q^{k-n/2}}
\left\{P_C(T)+q^{n/2+1-d}P_C\left(\frac{1}{qT}\right)T^{n+2-2d}\right\}.
\end{equation}
It satisfies $\deg \tilde P_C =2 \tilde g$ and the functional equation
\begin{equation}\label{eq:func_eq_tilde}
\tilde P_C(T)=\tilde P_C\left(\frac{1}{qT}\right)q^{\tilde g}T^{2\tilde g}
\end{equation}
where $\tilde g:= n/2-1-\min(d, d^\perp)$. 
\end{thm}
\prf By Definition \ref{dfn:zeta}, We have
\begin{equation}\label{eq:zeta_C}
\frac{P_C(T)}{(1-T)(1-qT)}(y(1-T)+xT)^n=\cdots +\frac{W_C(x,y)-x^n}{q-1}T^{n-d}+ \cdots
\end{equation}
and 
\begin{equation}\label{eq:zeta_C_perp}
\frac{P_{C^\perp}(T)}{(1-T)(1-qT)}(y(1-T)+xT)^n=\cdots +\frac{W_{C^\perp}(x,y)-x^n}{q-1}T^{n-d^\perp}+ \cdots.
\end{equation}
We suppose $d\leq d^\perp$. Then (\ref{eq:zeta_C_perp}) multiplied by $q^{k-n/2}T^{d^\perp-d}$ becomes
\begin{equation}\label{eq:zeta_C_perp-2}
\frac{q^{k-n/2} P_{C^\perp}(T) T^{d^\perp-d}}{(1-T)(1-qT)}(y(1-T)+xT)^n
=\cdots +\frac{q^{k-n/2}(W_{C^\perp}(x,y)-x^n)}{q-1}T^{n-d}+ \cdots.
\end{equation}
We add (\ref{eq:zeta_C}) and (\ref{eq:zeta_C_perp-2}), then divide it by $1+q^{k-n/2}$. It gives
$$\frac{\{P_{C}(T)+ q^{k-n/2} P_{C^\perp}(T) T^{d^\perp-d}\}/(1+q^{k-n/2})}{(1-T)(1-qT)}(y(1-T)+xT)^n=\cdots +\frac{\tilde W_{C}(x,y)-x^n}{q-1}T^{n-d}+ \cdots.$$
Thus we have 
\begin{equation}\label{eq:tilde_P_C-0}
\tilde P_C(T)=\frac{1}{1+q^{k-n/2}}\left\{P_{C}(T)+ q^{k-n/2} P_{C^\perp}(T) T^{d^\perp-d}\right\}
\end{equation}
by the existence and uniqueness of the zeta polynomial. The polynomial $P_{C^\perp}(T)$ can be substituted by
$$P_{C^\perp}(T)=P_C\left(\frac{1}{qT}\right)q^g T^{g+g^\perp}$$
where
\begin{eqnarray}
g&=&n+1-k-d,\\
g^\perp&=&k+1-d^\perp.
\end{eqnarray}
These formulas come from the original Duursma theory (see Duursma \cite[p.59]{Du2}). Hence we have
\begin{equation}\label{eq:tilde_P_C-1}
\tilde P_C(T)=\frac{1}{1+q^{k-n/2}}
\left\{P_C(T)+q^{n/2+1-d}P_C\left(\frac{1}{qT}\right)T^{n+2-2d}\right\}.
\end{equation}
When $d \geq d^\perp$, similarly we have
\begin{equation}\label{eq:tilde_P_C-2}
\tilde P_C(T)=\frac{T^{d-d^\perp}}{1+q^{k-n/2}}
\left\{P_C(T)+q^{n/2+1-d}P_C\left(\frac{1}{qT}\right)T^{n+2-2d}\right\}.
\end{equation}
These two formulas give (\ref{eq:tilde_P}). The functional equation (\ref{eq:func_eq_tilde}) is obtained in a similar manner to that of Duursma \cite[p.119]{Du3}. As to $\deg \tilde P_C$, first we note that 
\begin{equation}\label{eq:deg_P_C}
\deg P_C = \deg P_{C^\perp}=g+g^\perp=n+2-d-d^\perp
\end{equation}
(see \cite[p.59]{Du2}). By (\ref{eq:tilde_P_C-0}), we have $\deg \tilde P_C=n+2-2d=2\tilde g$ when $d\leq d^\perp$. The case $d\geq d^\perp$ is similar. \qed

\medskip
\noindent\rem When $C^\perp=C$, we can easily verify that $\tilde P_C(T)=P_C(T)$. Thus we have extended Duursma's theory in such a way that the zeta functions for codes which are not self-dual have the functional equation. 
\section{The MDS codes}

We consider the case where $C$ is an MDS code in the construction of $\tilde W_C(x,y)$ in Section 2 and prove Theorem \ref{thm:RH_MDS}. An $[n,k,d]$ code $C$ is called an MDS (maximal distance separable) code if $d=n-k+1$ is satisfied, i.e., the equality holds in the Singleton bound $d\leq n-k+1$. If $C$ is MDS, then so is $C^\perp$ and it has the parameter $[n,n-k,n+2-d]$. The weight enumerator $W_C(x,y)$ of an MDS code $C$ is determined only by $n$, $d$ and $q$. It can be explicitly given in terms of binomial coefficients: 
\begin{thm}\label{thm:MDSWE}
Let $W_C(x,y)=\sum_{i=0}^n A_i x^{n-i}y^i$ be the weight enumerator of an  $[n,k,d=n-k+1]$ MDS code $C$. Then we have
$$A_i={n\choose i} \sum_{j=0}^{i-d}(-1)^j {i\choose j} (q^{i-d+1-j}-1).
\quad (i\geq d)$$
\end{thm}
\prf MacWilliams-Sloane \cite[p.320, Theorem 6]{MaSl}. \qed

\medskip
\noindent We allow $A_i$ to be negative and $q$ be arbitrary integer greater than one. Even in the case $W_C(x,y)$ does not represent the weight distribution of an actual code, we are interested in the polynomial itself and often call it an ``MDS polynomial''. {F}rom now on we assume $d, d^\perp \geq 2$. What is crucial for our discussion is the following: 
\begin{thm}\label{thm:MDSzeta}
Let $C$ be MDS. Then we have $P_C(T)=1$. 
\end{thm}
\prf See Duusrma \cite[Proposition 1]{Du2}. In cases where $q$ is not a prime power, it is straightforward. \qed

\medskip
\noindent Now we determine the range of $d$ and $n$. Both $C$ and $C^\perp$ are MDS, so we can assume $d\leq d^\perp$ without loss of generality. Since $d^\perp=n+2-d$, $d\leq d^\perp$ is equivalent to
\begin{equation}\label{eq:d_and_d_perp}
d\leq \frac{n}{2}+1.
\end{equation}
If equality holds in (\ref{eq:d_and_d_perp}), then $\tilde g=0$ and $\tilde P_C(T)$ is a constant ($\tilde W_C(x,y)$ is an MDS polynomial in the ring ${\bf C}[x,y]^{G_q}$ in this case. It can happen when $n$ is even). We exclude this case and have $d\leq (n+1)/2$. Duursma's theory requires $d, d^\perp \geq 2$, therefore $d$ and $n$ can assume the values with
\begin{equation}\label{eq:range_dn}
2\leq d \leq \frac{n+1}{2}.
\end{equation}
Now let $C$ be an (actual or virtual) MDS code. Then we have $P_C(T)=P_{C^\perp}(T)=1$ by Theorem \ref{thm:MDSzeta}. The zeta polynomial $\tilde P_C(T)$ of the invariant polynomial $\tilde W_C(x,y)$ is given by Theorem \ref{thm:tilde_P} as
\begin{equation}\label{eq:tilde_P_C_MDS}
\tilde P_C(T)=\frac{1}{1+q^{k-n/2}}(1+q^{n/2+1-d}T^{n+2-2d}).
\end{equation}
We can easily see that all the roots of (\ref{eq:tilde_P_C_MDS}) lie on the circle $|T|=1/\sqrt{q}$. {F}rom Theorem \ref{thm:MDSWE} and (\ref{eq:range_dn}), $\tilde W_C(x,y)$ is of the form $x^n+A_d x^{n-d}y^d+\cdots$ and $A_d\ne 0$. This completes the proof of Theorem \ref{thm:RH_MDS}. 
\section{The general Hamming codes}

For $r\geq 2$ and a prime power $q$, the general Hamming $[(q^r-1)/(q-1)=n, n-r, 3]$ code ${\rm Ham}(r,q)$ over $\ff q$ is the dual code of an $[n, r, q^{r-1}]$ simplex code over $\ff q$ (see Pless et al. \cite[p.316]{Pl2}). Therefore we have
\begin{eqnarray}\label{eq:WEsimplex}
W_{{\rm Ham}(r,q)^\perp}(x,y) &=& 
 x^n+(q-1)nx^{\frac{n-1}{q}}y^{\frac{(q-1)n+1}{q}}\\
&=& x^n+(q^r-1)x^{n-q^{r-1}}y^{q^{r-1}}.\nonumber
\end{eqnarray}
In this section we assume $r\geq 3$, allow $q$ to be any integer with $q\geq 2$ and determine explicitly the zeta polynomial $\tilde P_{r,q}(T) := \tilde P_{{\rm Ham}(r,q)}(T)$ of the invariant polynomial $\tilde W_{{\rm Ham}(r,q)}(x,y)\in {\bf C}[x,y]^{G_q}$ constructed in the manner of Section 2. For our purpose, it is easier to handle with $W_{{\rm Ham}(r,q)^\perp}(x,y)$ than $W_{{\rm Ham}(r,q)}(x,y)$, so we fix the notation as follows: 
$$\begin{array}{ll}
C={\rm Ham}(r,q)^\perp, &  C^\perp={\rm Ham}(r,q), \\[3mm]
\displaystyle n=\frac{q^r-1}{q-1} & (\mbox{the length of $C$ and $C^\perp$}),\\[3mm]
d=q^{r-1} & (\mbox{the minimum distance of ${\rm Ham}(r,q)^\perp$}).
\end{array}$$
First we deduce the zeta polynomial $P_C(T)=P_{{\rm Ham}(r,q)^\perp}(T)$. We use the notion of the normalized weight enumerator (see Duursma \cite[Definition 2]{Du2}): 
\begin{dfn}\label{dfn:norm_WE}
For a weight enumerator $A(x,y)$ of the form (\ref{eq:homogen}), the normalized weight enumerator $a(t)$ is defined by 
$$a(t)=\frac{1}{q-1}\sum_{i=d}^n A_i \Bigm/ {n \choose i} t^{i-d}.$$
\end{dfn}
The following theorem gives the relation between $A(x,y)$ and its zeta polynomial $P(T)$: 
\begin{thm}[Duursma]\label{thm:NWE_and_P}
The weight enumerator $A(x,y)$, its zeta polynomial $P(T)$ and the normalized weight enumerator $a(t)$ are related by 
$$\frac{P(T)}{(1-T)(1-qT)}(1-T)^{d+1}\equiv a\left( \frac{T}{1-T} \right) \pmod{T^{n-d+1}}.$$
\end{thm}
\prf See \cite[Theorem 2]{Du2}. \qed

\medskip
\noindent For our code $C={\rm Ham}(r,q)^\perp$, the normalized weight enumerator $a(t)$ is quite simple: 
\begin{lem}\label{lem:NWE_Ham_perp}
Let $a_{r,q}(t)$ be the normalized weight enumerator of ${\rm Ham}(r,q)^\perp$. Then 
$$a_{r,q}(t)=n \Bigm/ {n \choose {q^{r-1}}},$$
i.e., $a_{r,q}(t)$ is a constant. 
\end{lem}
The proof is easy from (\ref{eq:WEsimplex}). Using this lemma and Theorem \ref{thm:NWE_and_P}, we can deduce the explicit form of  $P_C(T)$: 
\begin{prop}\label{prop:zeta_simplex}
For $r\geq 3$ and $q\geq 2$, the zeta polynomial $P_C(T)=P_{{\rm Ham}(r,q)^\perp}(T)$ is given by
\begin{equation}\label{eq:zeta_simplex}
P_C(T)=N_{r,q} \left[ 1+\sum_{j=1}^{n-d-1} 
\left\{ {{j+d-1}\choose{d-1}} -q {{j+d-2}\choose{d-1}} \right\} T^j\right],
\end{equation}
where $n/{n \choose {q^{r-1}}}$. 
\end{prop}
\prf Lemma \ref{lem:NWE_Ham_perp} and Theorem \ref{thm:NWE_and_P} gives
\begin{eqnarray}
P_C(T) &\equiv& a\left( \frac{T}{1-T} \right) 
\frac{(1-T)(1-qT)}{(1-T)^{d+1}}\pmod{T^{n-d+1}}\nonumber\\
   &\equiv& N_{r,q} \frac{1-qT}{(1-T)^d}\pmod{T^{n-d+1}}.\label{eq:equiv_zeta_simplex}
\end{eqnarray}
We have
$$\deg P_C = n+2-d-3 = n-d-1 < n-d+1$$
(see (\ref{eq:deg_P_C})), so $P_C(T)$ coincides with the power series expansion of $N_{r,q}(1-qT)/(1-T)^d$ up to the term of $T^{n-d-1}$. By the  expansion $(1-T)^{-d}=\sum_{j=0}^\infty {{j+d-1}\choose{d-1}} T^j$, we have
\begin{equation}\label{eq:expansion_(1-qT)/(1-T)_d}
\frac{1-qT}{(1-T)^d}=1+\sum_{j=1}^\infty 
\left\{ {{j+d-1}\choose{d-1}} -q {{j+d-2}\choose{d-1}} \right\} T^j.
\end{equation}
This formula gives the desired result. \qed

\medskip
\noindent\rem In the formula (\ref{eq:expansion_(1-qT)/(1-T)_d}), ${{j+d-1}\choose{d-1}} -q {{j+d-2}\choose{d-1}}=0$ holds if and only if $j=n-d$. Thus the term of $T^{n-d}$ really vanishes in (\ref{eq:expansion_(1-qT)/(1-T)_d}). 

\medskip
The main theorem in this section is the following: 
\begin{thm}\label{thm:tilde_P_Ham}
For $r\geq 3$ and $q\geq 2$, the zeta polynomial $\tilde P_{r,q}(T) := \tilde P_{{\rm Ham}(r,q)}(T)$ is given by 
$$\tilde P_{r,q}(T) =\frac{N_{r,q}}{1+q^{r-n/2}}(F_1(T)-qF_2(T)),$$
where
\begin{eqnarray*}
F_1(T)&=&\sum_{i=0}^{n-d-1} {{n-i-2}\choose{d-1}}q^{i+2-n/2}T^i
+\sum_{i=d-3}^{n-4} {{i+2}\choose{d-1}} T^i,\\
F_2(T)&=&\sum_{i=0}^{n-d-2} {{n-i-3}\choose{d-1}}q^{i+2-n/2}T^i
+\sum_{i=d-2}^{n-4} {{i+1}\choose{d-1}} T^i.
\end{eqnarray*}
\end{thm}
\rem If $r=2$, both ${\rm Ham}(r,q)$ and ${\rm Ham}(r,q)^\perp$ are MDS codes and are treated in Section 3. 

\medskip
\noindent\prf Since $d=q^{r-1}\geq 3$ if $r\geq 3$ and $q\geq 2$, we have from Theorem \ref{thm:tilde_P}, 
\begin{equation}\label{eq:tilde_P_Ham-0}
\tilde P_{r,q}(T)=\frac{T^{d-3}}{1+q^{r-n/2}}
\left\{P_C(T)+q^{n/2+1-d}P_C\left(\frac{1}{qT}\right)T^{n+2-2d}\right\}.
\end{equation}
The remaining task is to describe each term in (\ref{eq:tilde_P_Ham-0}) explicitly. We have from Proposition \ref{prop:zeta_simplex}, 
\begin{eqnarray}
\frac{T^{d-3}}{1+q^{r-n/2}} P_C(T) &=& \frac{N_{r,q}}{1+q^{r-n/2}}
\left[ T^{d-3}+ \sum_{j=1}^{n-d-1} 
\left\{ {{j+d-1}\choose{d-1}} -q {{j+d-2}\choose{d-1}} \right\} T^{d+j-3}\right]\nonumber\\
   &=& \frac{N_{r,q}}{1+q^{r-n/2}}
\left[ T^{d-3}+ \sum_{i=d-2}^{n-4} 
\left\{ {{i+2}\choose{d-1}} -q {{i+1}\choose{d-1}} \right\} T^i\right]
\label{eq:former_tilde_P_Ham}
\end{eqnarray}
by putting $d+j-3=i$. Next we have from Proposition \ref{prop:zeta_simplex} again that

\medskip
$\displaystyle \frac{T^{d-3}}{1+q^{r-n/2}} \cdot q^{n/2+1-d}P_C\left(\frac{1}{qT}\right)T^{n+2-2d}$
\begin{eqnarray}
 &=& \frac{N_{r,q} \, q^{n/2+1-d}}{1+q^{r-n/2}}
\left[ 1+ \sum_{j=1}^{n-d-1} 
\left\{ {{j+d-1}\choose{d-1}} -q {{j+d-2}\choose{d-1}} \right\} q^{-j}T^{-j}\right] T^{n-d-1}\nonumber\\
   &=& \frac{N_{r,q}}{1+q^{r-n/2}}
\left[ q^{n/2+1-d}T^{n-d-1} \phantom{\sum_{j=1}^{n-d-1}}\right.\nonumber\\
 & & + \left. \sum_{j=1}^{n-d-1} 
\left\{ {{j+d-1}\choose{d-1}} -q {{j+d-2}\choose{d-1}} \right\} q^{n/2+1-d-j}T^{n-d-j-1}\right].
\label{eq:latter_tilde_P_Ham}
\end{eqnarray}
By substitution $n-d-j-1=i$, (\ref{eq:latter_tilde_P_Ham}) equals
\begin{equation}\label{eq:latter_tilde_P_Ham-2}
\frac{N_{r,q}}{1+q^{r-n/2}} 
\left[ \sum_{i=0}^{n-d-2} 
\left\{ {{n-i-2}\choose{d-1}} -q {{n-i-3}\choose{d-1}} \right\}
q^{i+2-n/2}T^{i}+q^{n/2+1-d}T^{n-d-1}\right].
\end{equation}
The formulas (\ref{eq:tilde_P_Ham-0}), (\ref{eq:former_tilde_P_Ham}) and (\ref{eq:latter_tilde_P_Ham-2}) give
\begin{eqnarray}
\tilde P_{r,q}(T)&=&\frac{N_{r,q}}{1+q^{r-n/2}} 
\left[ \sum_{i=0}^{n-d-2} 
\left\{ {{n-i-2}\choose{d-1}} -q {{n-i-3}\choose{d-1}} \right\}
q^{i+2-n/2}T^{i}\right.\nonumber\\
   &+& \left. q^{n/2+1-d}T^{n-d-1} + T^{d-3}+
\sum_{i=d-2}^{n-4} 
\left\{ {{i+2}\choose{d-1}} -q {{i+1}\choose{d-1}} \right\} T^i\right].
\label{eq:tilde_P_Ham-1}
\end{eqnarray}
We make $F_1(T)$ by gathering positive terms in (\ref{eq:tilde_P_Ham-1}) and $F_2(T)$ from negative ones. \qed

\medskip
Theorem \ref{thm:RH_hamming} claims that all the roots of $\tilde P_{r,q}(T)$ above lie on the circle $|T|=1/\sqrt q$ if $q\geq4$. This is proved in several steps. We consider ``normalized'' zeta polynomial $\tilde P_{r,q}(T/\sqrt{q})$. Then the Riemann hypothesis is equivalent to the fact that all the roots of 
$\tilde P_{r,q}(T/\sqrt{q})$ lie on the unit circle. On the other hand, $\tilde P_{r,q}(T/\sqrt{q})$ is self-reciprocal, which is the result of the functional equation (\ref{eq:func_eq}) ($\sum_{i=0}^\nu a_i T^i$ is called self-reciprocal if $a_i=a_{\nu-i}$ for all $i$). Moreover, if $q\geq4$, $\tilde P_{r,q}(T/\sqrt{q})$ turns out to be of the form 
$$\tilde P_{r,q}(T/\sqrt{q})=
a_0+a_1 T+ \cdots +a_k T^k + a_k T^{m-k}+a_{k-1} T^{m-k+1} + \cdots +a_0 T^m$$
with $m>2k$ and $a_0>a_1> \cdots > a_k>0$. We can prove that all the roots of a self-reciprocal polynomial of this form lie on the unit circle using several results of classical function theory (an analogue of the Enestr\"om-Kakeya theorem, see Theorem \ref{thm:enest_kake}). We state the proof in the next two sections. 

\medskip
\noindent\rem We can also prove directly that $F_1(T/\sqrt q)$, $F_2(T/\sqrt q)$ and $\tilde P_{r,q}(T/\sqrt{q})$ are self-reciprocal, using the expressions in Theorem \ref{thm:tilde_P_Ham}. 
\section{An analogue of the Enestr\"om-Kakeya theorem}

In this section, we prove Theorem \ref{thm:ext_enest_kake}. This is a self-reciprocal analogue of the following theorem and our proof of Theorem \ref{thm:ext_enest_kake} is based on it: 
\begin{thm}[Enestr\"om-Kakeya]\label{thm:enest_kake}
Let $f(T)=a_0+a_1 T+ \cdots +a_k T^k$ satisfy $a_0>a_1> \cdots > a_k>0$. Then $f(T)$ has no roots in $|T|\leq 1$. 
\end{thm}
\prf Marden \cite[p.151, Exercise 4]{Ma}. \qed

\medskip
Now, suppose a self-reciprocal polynomial
\begin{equation}\label{eq:self-reci}
f(T)=a_0+a_1 T+ \cdots +a_k T^k + a_k T^{m-k}+a_{k-1} T^{m-k+1} + \cdots +a_0 T^m\quad (m>2k)
\end{equation}
satisfies $a_0>a_1> \cdots > a_k>0$. We write $f(T)$ as a sum of two polynomials $P(T)$ and $Q(T)$: 
\begin{eqnarray}
P(T) &:=& a_k T^{m-k}+a_{k-1} T^{m-k+1}+ \cdots +a_0 T^m,\nonumber\\
Q(T) &:=& a_0 +a_1 T+ \cdots +a_k T^k,\label{eq:P(T)Q(T)}
\end{eqnarray}
so $f(T)=P(T)+Q(T)$. Then, by the assumption $a_0>a_1> \cdots > a_k>0$, we can see from Theorem \ref{thm:enest_kake} that $Q(T)$ has no roots in $|T|\leq 1$. We apply Rouch\'e's theorem to $f(T)$. For simplicity, we state it in a restricted form: 
\begin{thm}\label{thm:rouche}
Let $C$ be a circle in ${\bf C}$, $D$ be the inside of $C$. Suppose functions $P(T)$ and $Q(T)$ are holomorphic in $C \cup D$ and $|P(T)|<|Q(T)|$ on $C$. Then $Q(T)$ and $P(T)+Q(T)$ have the same number of zeros in $D$. 
\end{thm}
\prf Ahlfors \cite[p.153, Corollary]{Ah}. See also Lehmer \cite[Lemma 3]{Le}. \qed

\medskip
For our polynomials $P(T)$ and $Q(T)$, we can prove the following: 
\begin{thm}\label{thm:ourP_and_Q}
We have $|P(T)|<|Q(T)|$ on $|T|=r$ for any $r$ with $0<r<1$. 
\end{thm}
By this theorem, we can see that $Q(T)$ and $f(T)=P(T)+Q(T)$ have the same number of roots in $|T|<r$. By Theorem \ref{thm:enest_kake} again, $f(T)$ has no roots in $|T|<r$. Since $r$ is arbitrary in $0<r<1$, we can verify that $f(T)$ has no roots in $|T|<1$. Now recall that $f(T)$ is self-reciprocal. We have 
$$T^m f\left( \frac{1}{T} \right) =f(T).$$
{F}rom this formula, we see that there is a one-to-one correspondence between a root in $|T|<1$ and that in $|T|>1$. We can conclude that $f(T)$ has no roots also in $|T|>1$, and all the roots of $f(T)$ lie on $|T|=1$. Hence we get Theorem \ref{thm:ext_enest_kake}. 

\medskip
\noindent{\bf Proof of Theorem \ref{thm:ourP_and_Q}.}

\medskip
First we need the following: 
\begin{lem}[Lagrange's identity]\label{lem:lagrange_id}
For any $A_i, B_i \in {\bf C}$, we have
$$|\sum_{i=0}^k A_i B_i |^2=\sum_{i=0}^k |A_i|^2 \sum_{i=0}^k |B_i|^2 
-\sum_{0\leq i<j \leq k} |A_i \overline{B_j} - A_j \overline{B_i}|^2.$$
\end{lem}
\prf Ahlfors \cite[p.9, Exercise 5]{Ah}. \qed

\medskip
\noindent Using this, we can prove the following: 
\begin{lem}\label{lem:P(T)Q(T)on1}
For $P(T)$ and $Q(T)$ in (\ref{eq:P(T)Q(T)}), we have
$$|P(T)|=|Q(T)|$$
on $|T|=1$. 
\end{lem}
\prf By letting $A_i=a_i$ and $B_i=T^i$ in Lemma \ref{lem:lagrange_id}, we get
\begin{equation}\label{eq:Q(T)squared}
|Q(T)|^2=(k+1)(a_0^2+ \cdots +a_k^2)
-\sum_{0\leq i<j \leq k} |a_i-a_jT^{j-i}|^2
\end{equation}
since $|T|=1$. As to $P(T)$, noting that $|P(T)|=|a_k + a_{k-1}T + \cdots +a_0 T^k|$ on $|T|=1$, we have
\begin{equation}\label{eq:P(T)squared}
|P(T)|^2=(k+1)(a_0^2+ \cdots +a_k^2)
-\sum_{0\leq i<j \leq k} |a_{k-i}-a_{k-j}T^{j-i}|^2
\end{equation}
by letting $A_i=a_{k-i}$ and $B_i=T^i$ in Lemma \ref{lem:lagrange_id}. By change of suffices in the sum in (\ref{eq:P(T)squared}), we have
\begin{equation}\label{eq:P(T)sum}
\sum_{0\leq i<j \leq k} |a_{k-i}-a_{k-j}T^{j-i}|^2
=\sum_{0\leq j'<i' \leq k} |a_{i'}-a_{j'}T^{i'-j'}|^2
=\sum_{0\leq i<j \leq k} |a_j-a_iT^{j-i}|^2.
\end{equation}
We compare the term for $(i,j)$ in (\ref{eq:Q(T)squared}) and (\ref{eq:P(T)sum}): 
\begin{eqnarray*}
|a_i-a_jT^{j-i}|^2-|a_j-a_iT^{j-i}|^2 &=& 
(a_i-a_jT^{j-i})(a_i-a_j \overline{T}^{j-i})-
(a_j-a_iT^{j-i})(a_j-a_i \overline{T}^{j-i})\\
&=&0
\end{eqnarray*}
since $|T|=1$. We see that the sums in the right hand sides in (\ref{eq:Q(T)squared}) and (\ref{eq:P(T)squared}) are the same, and we obtain $|P(T)|=|Q(T)|$ on $|T|=1$. \qed

\medskip
\noindent The proof of Theorem \ref{thm:ourP_and_Q} is completed by invoking the following well-known result:
\begin{thm}[The maximum principle]\label{thm:maxvalue}
Let $g(T)$ be holomorphic and nonconstant in a bounded (open) region $D\subset {\bf C}$ and continuous in $\overline{D}$ (the closure of $D$). Then $|g(T)|$ has its maximum $M$ on $\overline{D}-D$ and we have
$$|g(T)|<M$$
in $D$. 
\end{thm}
\prf Ahlfors \cite[p.134]{Ah}. \qed

\medskip
\noindent We apply Theorem \ref{thm:maxvalue} to $g(T):=P(T)/Q(T)$ and $D:=\{T \in {\bf C} \ ;\ |T|<1 \}$. Clearly $g(T)$ is meromorphic and nonconstant. It has no pole in $\overline{D}$ by Theorem \ref{thm:enest_kake}. Moreover, from Lemma \ref{lem:P(T)Q(T)on1}, $|g(T)|=1$ on the boundary of $D$. Therefore $|g(T)|<1$ in $D$ by Theorem \ref{thm:maxvalue} and we get Theorem \ref{thm:ourP_and_Q}. 
\section{Proof of Theorem \ref{thm:RH_hamming}}

In this section, we prove Theorem \ref{thm:RH_hamming}. We have from Theorem \ref{thm:tilde_P_Ham}, 
$$\tilde P_{r,q}\left(\frac{T}{\sqrt{q}}\right)
=\frac{N_{r,q}}{1+q^{r-n/2}}\left(F_1\left(\frac{T}{\sqrt{q}}\right)
-q F_2\left(\frac{T}{\sqrt{q}}\right) \right)$$
where 
\begin{eqnarray}
F_1\left(\frac{T}{\sqrt{q}}\right) &=&
\sum_{i=0}^{n-d-1} {{n-i-2}\choose{d-1}} q^{(i-n)/2+2} T^i
 + \sum_{i=d-3}^{n-4} {{i+2}\choose{d-1}} q^{-i/2} T^i,
\label{eq:F_1_self_recip}\\
F_2\left(\frac{T}{\sqrt{q}}\right) &=&
\sum_{i=0}^{n-d-2} {{n-i-3}\choose{d-1}} q^{(i-n)/2+2} T^i
 + \sum_{i=d-2}^{n-4} {{i+1}\choose{d-1}} q^{-i/2} T^i. 
\label{eq:F_2_self_recip}
\end{eqnarray}
Note that $n-d-1<d-3$ if $r\geq 3$ and $q\geq 4$. So there is no term of the same degree from two summations in $F_1(T)$ of Theorem \ref{thm:tilde_P_Ham}. Moreover, $\tilde P_{r,q}(T/\sqrt{q})$ is self-reciprocal. It follows from the functional equation (\ref{eq:func_eq}), but we can verify it directly by showing $F_1(T/\sqrt{q})$ and $F_2(T/\sqrt{q})$ are self-reciprocal. Hence we can assume $\tilde P_{r,q}(T/\sqrt{q})$ is of the form (\ref{eq:self-reci}). Let
\begin{equation}\label{eq:tilde_P_self_recip}
\frac{1+q^{r-n/2}}{N_{r,q}} \tilde P_{r,q}(T/\sqrt{q})=
a_0 + a_1 T + \cdots + a_{n-d-1}T^{n-d-1}
+a_{n-d-1} T^{d-3} + \cdots + a_0 T^{n-4}.
\end{equation}
\begin{lem}\label{lem:a_n-d-1}
If $r\geq 3$ and $q\geq 4$, we have 
$$a_{n-d-2} > a_{n-d-1} >0.$$
\end{lem}
\prf Recall $d=q^{r-1}$. Using the expressions (\ref{eq:F_1_self_recip}) and (\ref{eq:F_2_self_recip}), we have
$$a_{n-d-2}=q^{2-d/2}(q^{r-2}-1)$$
and 
$$a_{n-d-1}=q^{(3-d)/2}.$$
Therefore, $a_{n-d-1}>0$ and 
$$a_{n-d-2} - a_{n-d-1}=q^{(3-d)/2}\{ \sqrt{q}(q^{r-2}-1)-1\}.$$
The last expression is positive if $r\geq 3$ and $q\geq 4$. \qed
\begin{lem}\label{lem:a_i}
If $r\geq 3$ and $q\geq 4$, we have 
$$a_i > a_{i+1}$$
for $0\leq i \leq n-d-3$. 
\end{lem}
\prf Using the expressions (\ref{eq:F_1_self_recip}) and (\ref{eq:F_2_self_recip}), we have
\begin{eqnarray*}
a_i &=& {{n-i-2}\choose{d-1}}q^{(i-n)/2+2}-{{n-i-3}\choose{d-1}}q^{(i-n)/2+3},\\
a_{i+1} &=& {{n-i-3}\choose{d-1}}q^{(i-n+5)/2}-{{n-i-4}\choose{d-1}}q^{(i-n+7)/2}.
\end{eqnarray*}
{F}rom these formulas, we have
$$\left\{ q^{(i-n)/2+2}{{n-i-4}\choose{d-1}} \right\}^{-1}
(n-d-i-1)(n-d-i-2)(a_i-a_{i+1})$$
$$=(n-i-2)(n-i-3)-(n-i-3)(n-d-i-1)(q+\sqrt{q})
+(n-d-i-1)(n-d-i-2)q\sqrt{q}.$$
It suffices to show that the right hand side is positive. It is a quadratic function of the parameter $i$, so we denote it by $g(i)=ai^2+bi+c$. We can show that $a,b,c >0$ if $q\geq 4$. Indeed, first we have 
$$a=q\sqrt{q}+1-(q+\sqrt{q})=(\sqrt{q}-1)(q-1)>0$$
if $q\geq 2$. As to $b$, recall $n=(q^r-1)/(q-1)$ and $d=q^{r-1}$. We have 
$$\sqrt{q}(q-1)b=
q^r(\sqrt{q}-1)(q-1)+q^{3/2}(3 q^{3/2}-4q-5\sqrt{q}+7)+\sqrt{q}(2\sqrt{q}-3)$$
(such calculation can be easily done with the help of some expression manipulation program). As above, $q^r(\sqrt{q}-1)(q-1)>0$ if $q\geq 2$ and $2\sqrt{q}-3>0$ if $q\geq 3$. We can easily show $3 q^{3/2}-4q-5\sqrt{q}+7>0$ if $q\geq 4$ (show that $(3 q^{3/2}-4q-5\sqrt{q}+7)|_{q=4}>0$ and  $(3 q^{3/2}-4q-5\sqrt{q}+7)'>0$ in $q\geq 4$). Therefore $b>0$ if $q\geq 4$. 

We can similarly show $c>0$. Because
\begin{eqnarray*}
\sqrt{q}(q-1)^2c &=& q^r\{q(q^{3/2}-2q-2q^{1/2}+4)+(q^{1/2}-2)\}\\
   &+& q^{5/2}(2q^{3/2}-3q-4q^{1/2}+7)\\
   &+& q^{1/2}(q^2+2q^{3/2}-7q+2),
\end{eqnarray*}
and we can show that $q^{3/2}-2q-2q^{1/2}+4\geq 0$ if $q\geq 4$, and that all other functions in the parentheses ( ) are positive in $q\geq 4$. 

It follows that $g(i)>0$ if $i\geq 0$ and $a_i-a_{i+1}>0$. \qed

\medskip
We can conclude from Lemmas \ref{lem:a_n-d-1} and \ref{lem:a_i} that 
\begin{equation}\label{eq:coeff_tilde_P}
a_0>a_1> \cdots > a_{n-d-1}>0
\end{equation}
for the coefficients in (\ref{eq:tilde_P_self_recip}). The assumption of Theorem \ref{thm:ext_enest_kake} is satisfied and the proof of Theorem \ref{thm:RH_hamming} is completed. 

\medskip
\noindent\rem (1) If we write $\tilde W_{{\rm Ham}(r,q)}(x,y)$ in the form (\ref{eq:homogen}), $n=(q^r-1)/(q-1)$ and $d=3$. Thus we have found infinitely many invariant polynomials satisfying the Riemann hypothesis for a small $d$. 

\smallskip
\noindent (2) When $q=2,3$, the coefficients of $\tilde P_{r,q}(T/\sqrt{q})$ does not satisfy (\ref{eq:coeff_tilde_P}) as the following examples show. So, in these cases, we cannot prove the Riemann hypothesis in a method described so far, but numerical experiments imply that it is very plausible that the Riemann hypothesis is true also for $q=2,3$. 
\begin{exam}\label{exam:q=2}\rm 
(i) Let $r=3$, $q=2$. Then 
$$\tilde W_{{\rm Ham}(3,2)}(x,y)=
x^7+\frac{7}{1+\sqrt{2}}x^4y^3+7x^3y^4+\frac{7}{1+\sqrt{2}}y^7,$$
\begin{eqnarray*}
F_1(T)-2 F_2(T) &=& \frac{1}{\sqrt{2}} + (1+\sqrt{2})T
 + (2+\sqrt{2})T^2 + 2T^3,\\
F_1\left(\frac{T}{\sqrt{2}}\right)-2 F_2\left(\frac{T}{\sqrt{2}}\right)
 &=& \frac{1}{\sqrt{2}} + \left(1+\frac{1}{\sqrt{2}}\right)T
 + \left(1+\frac{1}{\sqrt{2}}\right)T^2
 + \frac{1}{\sqrt{2}}T^3\\
 &=& \frac{1}{\sqrt{2}}(T+1)\left(T-\frac{-1+i}{\sqrt{2}}\right)
 \left(T-\frac{-1-i}{\sqrt{2}}\right).
\end{eqnarray*}
Hence $\tilde W_{{\rm Ham}(3,2)}(x,y)$ satisfies the Riemann hypothesis. 

\medskip
\noindent(ii) Let $r=4$, $q=2$. Then $\tilde P_{4,2}(T/\sqrt{q})$ is of degree 11. We normalize $\tilde P_{4,2}(T/\sqrt{q})$ with a suitable constant $C$ as $C\tilde P_{4,2}(T/\sqrt{q})=1+a_1T+\cdots$. Then we can approximate the coefficients as follows:
$$\begin{array}{ll}
a_0=1                   & a_3 \approx 1.028518954 \\
a_1 \approx 1.414213562 & a_4 \approx 0.606060606 \\
a_2 \approx 1.363636363 & a_5 \approx 0.317735799\\
\end{array}$$
and $a_6, \cdots, a_{11}$ are the same as above in the reverse order. We have $a_0<a_1>a_2>a_3>a_4>a_5>0$, but according to the numerical experiment, $\tilde W_{{\rm Ham}(4,2)}(x,y)$ seems to satisfy the Riemann hypothesis. 

\medskip
\noindent(iii) Let $r=5$, $q=2$. Then $\tilde P_{5,2}(T/\sqrt{q})$ is of degree 27. Choose $C$ as $C\tilde P_{5,2}(T/\sqrt{q})=1+a_1T+\cdots$. Then the coefficients are approximated as
$$\begin{array}{lr}
a_0=1                   &    a_7 \approx 0.2623468638 \\
a_1 \approx 1.414213562 &    a_8 \approx 0.1391304348 \\
a_2 \approx 1.444444444 &    a_9 \approx 0.0655867159\\
a_3 \approx 1.257078722 & a_{10} \approx 0.0268497330\\
a_4 \approx 0.977777778 & a_{11} \approx 0.0092051531\\
a_5 \approx 0.691393297 & a_{12} \approx 0.0024887453\\
a_6 \approx 0.446376812 & a_{13} \approx 0.0005216551\\
\end{array}$$
and $a_{14}, \cdots, a_{27}$ are the same as above in the reverse order. In this case we have $a_0<a_1<a_2>a_3> \cdots >a_{13}>0$. The Riemann hypothesis seems to be true. 
\end{exam}
\begin{exam}\label{exam:q=3}\rm 
Let $r=3$, $q=3$. Then $\tilde P_{3,3}(T/\sqrt{q})$ is of degree 9. Choose $C$ as $C\tilde P_{3,3}(T/\sqrt{q})=1+a_1T+\cdots$. Then $a_1\approx 1.039230485$, $a_2=0.6$, $a_3\approx 0.1732050808$, $a_4=a_5=0$, and $a_6, \cdots, a_9$ are the same but in the reverse order. The Riemann hypothesis seems to be true. In many other ${\rm Ham}(r,3)$ with $r\geq 4$, we can observe that the coefficient of $T$ in $\tilde P_{r,3}(T/\sqrt{q})$ is greater than the constant term. 
\end{exam}
\section{The Golay codes}

In this section we consider the case where $C={\cal G}_{23}$, the binary $[23,12,7]$ Golay code or $C={\cal G}_{11}$, the ternary $[11,6,5]$ Golay code. We have 
\begin{eqnarray}
W_{{\cal G}_{23}}(x,y)&=&x^{23}+253x^{16}y^7+506x^{15}y^8+1288x^{12}y^{11}+1288x^{11}y^{12}+506x^{8}y^{15}\nonumber\\
  & &+253x^7y^{16}+y^{23},\label{eq:WE_g23}\\
W_{{\cal G}_{11}}(x,y)&=&x^{11}+132x^6y^5+132x^5y^6+330x^3y^8+110x^2y^9+24y^{11}
\label{eq:WE_g11}
\end{eqnarray}
(see \cite[p.94]{Pl2}) or
\begin{eqnarray}\label{eq:WE_g23_perp}
W_{{{\cal G}_{23}}^\perp}(x,y)&=&x^{23}+506x^{15}y^8+1288x^{11}y^{12}+253x^7y^{16},\\
W_{{{\cal G}_{11}}^\perp}(x,y)&=&x^{11}+132x^5y^6+110x^2y^9.
\label{eq:WE_g11_perp}
\end{eqnarray}
The case $C={\cal G}_{11}$ is quite easy: 
\begin{prop}\label{prop:g11}
The zeta polynomial $\tilde P_{{\cal G}_{11}}(T)$ of the invariant polynomial $\tilde W_{{\cal G}_{11}}(x,y)$ is given by 
$$\tilde P_{{\cal G}_{11}}(T)=\frac{\sqrt{3}-1}{14}(\sqrt{3}T+1)(3T^2+3T+1).$$
All the roots of $\tilde P_{{\cal G}_{11}}(T)$ lie on the circle $|T|=1/\sqrt{3}$. 
\end{prop}
\prf The explicit form of $\tilde P_{{\cal G}_{11}}(T)$ can be obtained by computer calculation. The latter statement is obvious. \qed

\medskip
Next we consider $C={\cal G}_{23}$. By computer calculation we get the following: 
\begin{prop}\label{prop:g23_tilde_P}
Let $\tilde P_{{\cal G}_{23}}(T)$ be the zeta polynomial of the invariant polynomial $\tilde W_{{\cal G}_{23}}(x,y)$. Then 
\begin{eqnarray}
\frac{25194}{2-\sqrt{2}}\tilde P_{{\cal G}_{23}}\left(\frac{T}{\sqrt{2}}\right)
  &=&13(1+T^{11})+(13+39\sqrt{2})(T+T^{10})\nonumber\\
  &+&(130+39\sqrt{2})(T^2+T^9)+(130+156\sqrt{2})(T^3+T^8)\nonumber\\
  &+&\frac{591+312\sqrt{2}}{2}(T^4+T^7)+\frac{591+459\sqrt{2}}{2}(T^5+T^6).
\label{eq:g23_tilde_P}
\end{eqnarray}
\end{prop}
The polynomial $\tilde P_{{\cal G}_{23}}(T/\sqrt{2})$ does not satisfy the assumption of Theorem \ref{thm:ext_enest_kake}. We would like to verify the Riemann hypothesis for $\tilde W_{{\cal G}_{23}}(x,y)$ as theoretically as possible. Our method is influenced by that of Duursma \cite[Section 5]{Du3}
 and \cite[Section 5]{Du4}. 

Let $f(x)$ be a real polynomial of degree $n$. Then $f^\ast(T)=T^n f((T+T^{-1})/2)$ is a self-reciprocal polynomial of degree $2n$. If $T=e^{i\theta}$, then $(T+T^{-1})/2=\cos \theta$, so the behavior of $f^\ast(T)$ on the unit circle can be captured by the behavior of $f(x)$ in the interval $[-1,1]$. We denote this mapping by $\rho$: 
\begin{equation}\label{eq:rho_def}
\rho: f\mapsto f^\ast.
\end{equation}
We would like to pull back $f(x)$ from a given self-reciprocal polynomial $f^\ast(T)$. It turns out that the inverse mapping $\rho^{-1}$ always exists. To clarify this, we introduce two linear spaces of polynomials:
\begin{eqnarray*}
V_n &:=& \{ a_0 + a_1 x + \cdots + a_n x^n \ ;\ a_j\in{\bf R} \},\\
W_n &:=& \{b_0 + b_1 T + \cdots + b_n T^n + b_{n-1} T^{n+1} + \cdots b_0 T^{2n}\ ;\ b_j\in{\bf R} \}.
\end{eqnarray*}
The operations in the vector spaces are the same as ordinary summation of polynomials and multiplication by real numbers. 
\begin{lem}\label{lem:rho_isom}
The mapping $\rho: V_n \to W_n$ defined by (\ref{eq:rho_def}) is a linear isomorphism. 
\end{lem}
\prf Clearly $\rho$ is linear. For the later use, we describe the matrix $A_\rho$ of $\rho$ with respect to the bases 
\begin{eqnarray*}
V_n &=& {\bf R}[1, x, \cdots, x^n],\\
W_n &=& {\bf R}[T^n, T^{n-1}+T^{n+1}, \cdots 1+T^{2n}].
\end{eqnarray*}
Because $\rho$ maps $1, x, x^2, \cdots$ as
$$\begin{array}{lcl}
1 &\mapsto& T^n,\\
x &\mapsto& \frac{{1\choose 0}}{2}(T^{n-1}+T^{n+1}),\\
x^2 &\mapsto& 
  \frac{{2\choose 0}}{2^2}(T^{n-2}+T^{n+2})+\frac{{2\choose 1}}{2^2}T^n,\\
x^3 &\mapsto& 
  \frac{{3\choose 0}}{2^3}(T^{n-3}+T^{n+3})
  +\frac{{3\choose 1}}{2^3}(T^{n-1}+T^{n+1}),\\
x^4 &\mapsto& 
  \frac{{4\choose 0}}{2^4}(T^{n-4}+T^{n+4})
  +\frac{{4\choose 1}}{2^4}(T^{n-2}+T^{n+2})
  +\frac{{4\choose 2}}{2^4}T^n,\\
\multicolumn{3}{c}{\dotfill}.
\end{array}$$
we have for even $n$, 
$$A_\rho=\left[\begin{array}{ccccccc}
1 & 0 & 2^{-2}{2\choose 1} & 0 & 2^{-4}{4\choose 2} & \cdots & 2^{-n}{n\choose{n/2}}\\
 & 2^{-1}{1\choose 0} & 0 & 2^{-3}{3\choose 1} & 0 & \cdots & 0\\
 & & 2^{-2}{2\choose 0} & 0 & 2^{-4}{4\choose 1} & \cdots & 2^{-n}{n\choose{n/2-1}}\\
 & & & 2^{-3}{3\choose 0} & 0 & \cdots & 0\\
 & & & & 2^{-4}{4\choose 0} & \cdots & 2^{-n}{n\choose{n/2-2}}\\
 & & & & & \ddots & \vdots \\
 & & & & & & 2^{-n}{n\choose{0}}
\end{array}\right],$$
where all the elements in the lower triangular part are zero. If $n$ is odd, the $(n+1)$-th column is replaced by $^t[0, 2^{-n}{n\choose{(n-1)/2}}, 0,  2^{-n}{n\choose{(n-1)/2-1}}, \cdots, 2^{-n}{n\choose{0}}]$. By this expression, $\det A_\rho=2^{-n(n+1)/2}\ne 0$, so $A_\rho$ is regular. \qed

\medskip
We sketch the proof that all the roots of $\tilde P_{{\cal G}_{23}}(T/\sqrt{2})$ lie on the unit circle. The calculation is done with the help of a computer. Let $F(T)=(25194/(2-\sqrt{2}))\tilde P_{{\cal G}_{23}}(T^2/\sqrt{2})$. Then $\deg F=22$ and $F(T)$ is self-reciprocal. Construct $A_\rho$ for $n=11$ and map $F(T)$ by $\rho^{-1}$ using ${A_\rho}^{-1}$. Then we have
\begin{eqnarray*}
(\rho^{-1}F)(x) &=& 26624 x^{11}+(-66560+19968\sqrt{2})x^9\\
  &+& (74880-39936\sqrt{2})x^7+(-45760+29952\sqrt{2})x^5\\
  &+& (14324-9984\sqrt{2})x^3+(-1754+1239\sqrt{2})x.
\end{eqnarray*}
We can also verify $(\rho^{-1}F)(-1)<0$, $(\rho^{-1}F)(-0.8)>0$, $(\rho^{-1}F)(-0.6)<0$, $(\rho^{-1}F)(-0.4)>0$, $(\rho^{-1}F)(-0.2)<0$, $(\rho^{-1}F)(-0.1)>0$. It follows that $(\rho^{-1}F)(x)$ has at least five roots in the interval $(-1,0)$. It has the same number of roots in $(0,1)$ since it is an odd function of $x$, and $(\rho^{-1}F)(0)=0$. Thus all the roots of $(\rho^{-1}F)(x)$ are distinct and lie in $(-1,1)$. Let $T=e^{i\theta}$. Then $x=\cos\theta$. While $x$ moves from $1$ to $-1$, $T^2$ goes once around the unit circle. Therefore $\tilde P_{{\cal G}_{23}}(T^2/\sqrt{2})$ assumes zero exactly 11 times on $|T|=1$. 
\section{Appendix --- an elementary proof of existence of $P(T)$}

Existence and uniqueness of the zeta polynomial $P(T)$ for a linear code was first established in \cite[Section 9]{Du1}, but a detailed proof is not given. Here we give an alternative, elementary proof, including the case $W(x,y)\in{\bf C}[x,y]$. 

Suppose $W(x,y)\in{\bf C}[x,y]$ is a polynomial of the form (\ref{eq:homogen}). First note that 

\bigskip
$\displaystyle f(T):=\frac{1}{(1-T)(1-qT)}(y(1-T)+xT)^n$
\begin{eqnarray*}
&=&(1+T+T^2+\cdots)(1+qT+q^2T^2+\cdots)((x-y)T+y)^n\\
&=&(1+c_1 T+c_2 T^2+\cdots)\Bigl\{\sum_{j=0}^n \Bigl({n \atop j}\Bigr) (x-y)^jy^{n-j}T^j\Bigr\}
\end{eqnarray*}
for some $c_j\in{\bf N}$. Expanding the last formula, we find for some integers $b_{ij}$, 
$$\begin{array}{lr}
\mbox{the constant term} & =y^n,\\
\mbox{the coefficient of } T& =nxy^{n-1}+(c_1-n)y^n,\\
\cdots\cdots & \cdots\cdots\cdots\\
\mbox{the coefficient of } T^i & =b_{i0}x^iy^{n-i}+b_{i1}x^{i-1}y^{n-i+1}+\cdots +b_{ii}y^n,\\
\cdots\cdots & \cdots\cdots\cdots\\
\mbox{the coefficient of } T^{n-d} & =b_{n-d,0}x^{n-d}y^{d}+b_{n-d,1}x^{n-d-1}y^{d+1}+\cdots +b_{n-d,n-d}y^n.
\end{array}$$
Let $a_0$, $a_1$, $\cdots$, $a_{n-d}\in{\bf C}$ and we form a function $F(T):=(a_0+a_1 T+\cdots +a_{n-d}T^{n-d})f(T)$. Then the coefficient of $T^{n-d}$ of $F(T)$ is 
\begin{eqnarray}
a_{n-d}y^n& &\nonumber\\
+a_{n-d-1}\{nxy^{n-1}+(c_1-n)y^n\}& &\nonumber\\
\cdots\cdots\cdots\qquad\qquad & &\nonumber\\
+a_i\{b_{i0}x^iy^{n-i}+b_{i1}x^{i-1}y^{n-i+1}+\cdots +b_{ii}y^n\}& &\nonumber\\
\cdots\cdots\cdots\qquad\qquad & &\nonumber\\
+a_0\{b_{n-d,0}x^{n-d}y^{d}+b_{n-d,1}x^{n-d-1}y^{d+1}+\cdots +b_{n-d,n-d}y^n\}.& &
\label{eq:coeff_t^n-d}
\end{eqnarray}
On the other hand, since $(W(x,y)-x^n)/(q-1)=(A_dx^{n-d}y^d+\cdots +A_ny^n)/(q-1)$, we can determine $a_0$, $a_1$, $\cdots$, $a_{n-d}$ so that (\ref{eq:coeff_t^n-d}) coincides with $(W(x,y)-x^n)/(q-1)$ (the system of linear equations for determining $a_0$, $a_1$, $\cdots$, $a_{n-d}$ has a regular coefficient matrix). So we can always determine the zeta polynomial $P(T)$ from a given $W(x,y)$ uniquely as $P(T)=a_0+a_1 T+\cdots +a_{n-d}T^{n-d}$. \qed

\end{document}